\documentclass[11pt]{amsart}
\usepackage{geometry}                
\usepackage{graphicx}
\usepackage{amssymb}
\usepackage{epstopdf}
 \newtheorem{theorem}{Theorem}[section]

     \newtheorem{conjecture}[theorem]{Conjecture}  
    \theoremstyle{definition}

\newtheorem{examples}[theorem]{Examples}

  \newcommand{\G}{\Gamma}

  \newcommand{\AG}{A_\Gamma}
  
   \newcommand{\OF}{\ensuremath{{\mathcal{O}(F_n)}}}
  \newcommand{\Z}{\ensuremath{{\mathbb{Z}}}}
\newcommand{\R}{\ensuremath{{\mathbb{R}}}}
\newcommand{\Hy}{\ensuremath{{\mathbb{H}}}}

  \newcommand{\C}{\ensuremath{{\mathbb{C}}}}
 \newcommand{\B}{\ensuremath{{\mathcal{B}}}}
  \newcommand{\Sa}{\ensuremath{{\mathcal{S}_\Gamma}}}

\begin{document}

\title[Right-angled Artin groups]
{An Introduction to Right-Angled Artin Groups}

\date{October 2006}
\author[R.~Charney]{Ruth Charney}
	\address{Department of Mathematics\\
	Brandeis University\\
	Waltham, MA 02454\\
	USA}
	\email{charney@brandeis.edu}
	\thanks{Charney was partially
      supported by NSF grant DMS-0405623.}

\keywords{ Artin group, CAT(0) cube complex}

\subjclass{20F36}

\begin{abstract}  Recently, right-angled Artin groups have attracted much attention in geometric group theory.  They have a rich structure of subgroups and nice algorithmic properties, and they give rise to cubical complexes with a variety of applications.  This survey article is meant to introduce readers to these groups and to give an overview of the relevant literature.

\end{abstract}

\maketitle

Artin groups span a wide range of groups from braid groups to free groups to free abelian groups and have strong connections with geometry.  They are defined by presentations of specific form and are closely related to Coxeter groups.  Right-angled Artin groups are those Artin groups for which all relators are commutators between specified generators. 
On first glance the most elementary class of Artin groups, right-angled Artin groups turn out to have a surprising richness and flexibility that has led to some remarkable applications.  

Right-angled Artin groups were first introduced  in the 1970's by A.~Baudisch \cite{Bau} and further developed in the 1980's by  C.~Droms \cite{Dro1} \cite{Dro2} \cite{Dro3} under the name ``graph groups".   They have been studied extensively since that time (as evidenced by the long bibliography to this article). Recently, these groups have attracted much interest in geometric group theory due to their actions on CAT(0) cube complexes.  

In this survey article, we introduce the reader to the basic algebraic and geometric properties of right-angled Artin groups.  While we have tried to touch on all the major developments and to give references to related work where appropriate, the author apologizes for any (no doubt many) omissions.
The first section of the paper gives a brief introduction to more general Artin groups.  The second section reviews essential properties of right-angled Artin (and Coxeter) groups and the cube complexes associated to them.  The third section introduces some areas of current research on automorphisms and quasi-isometries of right-angled Artin groups, and the final section discusses a few of our favorite applications.

This article grew out of a series of talks given at the Pacific Geometry Conference at Oregon State University in June 2006.  It is written for a broad audience.  We hope that the content is interesting to geometric group theorists, but still accessible to those outside the field.  

The author would like to thank John Meier for offering some excellent suggestions during the preparation of this article.


\section{A brief introduction to Artin groups}

In order to put the theory of right-angled Artin groups in context, we begin with a very brief introduction to more general Artin groups. Artin groups arose as a generalization of braid groups, so we begin our discussion with the classical braid groups. 

\subsection{Braid groups}

 Let $\B_n$ denote the braid group on $n$ strands.  An element of $\B_n$ can be represented by $n$ intertwining strings attached at the top and bottom to $n$ fixed points. Two such braids are considered equal if one can be isotoped to the other without moving the endpoints of the braid.  Multiplication in the braid group is given by concatenating braids, one below the other. 

The braid group $\B_n$ is generated by $n-1$ elements $\tau_1, \dots ,\tau_{n-1}$ where $\tau_i$ is the braid which crosses the $i^{th}$ string over the $(i+1){st}$ string.  This gives rise to a well-known presentation of the braid group
\[
\B_n=\langle \tau_1, \dots ,\tau_{n-1} \mid  \tau_i\tau_{i+1}\tau_i=\tau_{i+1}\tau_i\tau_{i+1}, \, \tau_i\tau_j=\tau_j\tau_i \, \textrm{for $i+1<j$}
 \rangle.
\]

\begin{figure}
\begin{center}
\includegraphics[height=2.5cm]{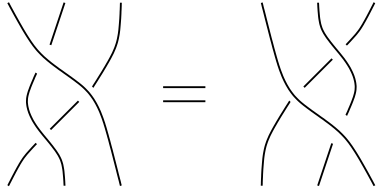}
\end{center}
\caption{The relation $\tau_i\tau_{i+1}\tau_i=\tau_{i+1}\tau_i\tau_{i+1}$ } 
\label{Fig1}
\end{figure}

If we add additional relations $\tau_i=\tau_i^{-1}$  to this presentation (so that over-crossings equal under-crossings), then the group reduces to the symmetric group $Sym_n$.  Thus, there is a natural surjection $\B_n \to Sym_n$.

The significance of the braid groups lie in their connections with geometry.  If we think of a braid as taking place inside a solid cylinder,  then taking horizontal slices of the cylinder gives a description of a braid as a set of $n$ points moving around in the disk, or equivalently, in the complex plane.  We can view this as a path in the configuration space $Y=\C^n - \mathbf{\Delta}$, where $\mathbf{\Delta}$ is the ``thick diagonal"
\[
\mathbf{\Delta} = \{ (z_1, \dots z_n)\in \C^n  \mid \textrm{$z_i=z_j$ for some $i\neq j$}\}.
\]
Up to permutation, the $n$ points begin and end in the same position, so this path forms a loop in $Y/ Sym_n$  where $Sym_n$ acts on $\C^n$ by permuting the coordinates.  It was shown by Fox and Neuwirth \cite{FoNe} that this correspondence gives and isomorphism
\[
\B_n \cong  \pi_1(Y/ Sym_n).
\]

The braid group is also isomorphic to the mapping class group of the $n$-punctured disk (where mapping classes are required to fix the boundary).   Intuitively this correspondence can be thought of as follows.  Given a homeomorphism of the $n$ punctured disk,  if we fill in the punctures, the homeomorphism can be isotoped back to the identity map, keeping the boundary fixed.  Watching the $n$ points move under this isotopy, we see a loop in the configuration space, or equivalently, a braid.

\subsection{Artin groups and the $K(\pi,1)$-conjecture}

Artin groups arose as natural generalizations of braid groups.
An \emph{Artin group} $A$ is a group with presentation of the form
\[
A=\langle s_1,\dots,s_n\ \mid\ 
\underbrace{s_is_js_i\dots}_{m_{ij}}=\underbrace{s_js_is_j\dots}_{m_{ij}}\ 
\textrm{ for all } i\neq j\rangle \,.
\]
where $m_{ij}=m_{ji}$ is an integer $\geq 2$ or $m_{ij}=\infty$ in which case we omit the relation between $s_i$ and $s_j$.  If we add to this presentation the additional relations $s_i=s_i^{-1}$ for all $i$, we obtain a \emph{Coxeter group}
\begin{align*}
W &= \langle  s_1,\dots,s_n \mid  s_i=s_i^{-1}, \, s_is_js_i\dots = s_js_is_j\dots \textrm{ for all } i\neq j \rangle \\
 &= \langle  s_1,\dots,s_n  \mid  (s_i)^{2}=1, \, (s_is_j)^{m_{ij}}=1\ \textrm{ for all } i\neq j\rangle .
\end{align*}
Let $S$ denote the generating set $\{s_1, \dots ,s_n\}$.  Then the pair $(A,S)$ is called an \emph{Artin system} and $(W,S)$ is called a \emph{Coxeter system}.  Every Artin system is associated to a Coxeter system and vice versa. By abuse of terminology, we will always assume that the generating set $S$ is specified when we refer to  a ``Coxeter group" or ``Artin group".

Any Coxeter group can be represented as a discrete reflection group, that is, a discrete group of linear transformations of a finite dimensional vector space $V$ with the generators $s_i$ acting as refections with respect to some bilinear form $B$.  (A ``reflection with respect to $B$" is a linear map $r_u(v)=v-2B(u,v) u$ for some unit vector $u$.)   If $W$ is a finite group, the associated bilinear form is positive definite, so we can identify $V,B$ with $\R^n$ equipped with the usual dot product.  In this case, we have a finite hyperplane arrangement in $\R^n$,
\[
\mathcal A = \{ H_r \mid \textrm{$H_r$ is the fixed set of some reflection $r \in W$} \}.
\]
The set of points in $\R^n$ with non-trivial isotropy under the $W$-action is precisely the union of these hyperplanes. Complexifying this action gives a finite arrangement of complex hyperplanes $\C H_r$ in $\C^n$ such that $W$ acts freely on the complement, 
$Y_W = \C^n - (\cup \, \C H_r) $.  These spaces arise in singularity theory and were studied in the 1970's by Brieskorn, Deligne, and others.  Brieskorn \cite{Bri}  first introduced Artin groups  when he  showed that the fundamental groups of the orbit spaces $Y_W/W$ have presentations of this form.   Shortly after that, Deligne \cite{Del} proved the following stronger theorem.

\begin{theorem}
Let $W$ be a finite Coxeter group and $A$ the associated Artin group.  Then $Y_W/W$ is aspherical with fundamental group $A$.  That is, $Y_W/W$ is a $K(A,1)$-space.
\end{theorem}

Based on work of Vinberg \cite{Vin},  one can define an analogous hyperplane complement $Y_W$ in $\C \otimes V$ for an infinite Coxeter group $W$.  (The definition is slightly different as one must restrict to the ``Tits cone,"  a $W$-equivariant open cone in $V$.)  The following conjecture is sometimes attributed to Arnold.

\begin{conjecture}[The $K(\pi,1)$-Conjecture]
The analogue of Deligne's theorem holds for infinite Coxeter groups.
\end{conjecture}

Part of this conjecture, namely that $A=\pi_1(Y_W/W)$, was proved for arbitrary Coxeter groups by van der Lek in \cite{vdL1}.  The conjecture that $Y_W/W$ is aspherical was proved for some, but not all, classes of infinite Coxeter groups by the author and M.~Davis in \cite{CD1} and \cite{Ch3}.  The conjecture remains open for arbitrary Coxeter groups.

In addition to their description as fundamental groups of hyperplanes complements, Artin groups associated to several of the finite and Euclidean Coxeter groups can be described as subgroups of mapping class groups or as ``orbifold braid groups".  See \cite{ChCr} and \cite{All} for a discussion of these groups. 

\medskip

We say that an Artin group is \emph{spherical} (or \emph{finite type}) if the associated Coxeter group is finite and \emph{non-spherical} (or \emph{infinite type}) otherwise.  Note that the Artin group itself is never finite; indeed, spherical Artin groups are torsion-free. (Non-spherical Artin groups are conjectured to be torsion-free.) Spherical Artin groups behave much like braid groups and we know a great deal about them.  By contrast, we know very little about non-spherical Artin groups except in some special cases (precisely the cases for which the $K(\pi,1)$ conjecture is known to hold). For example, spherical Artin groups are known to have the following properties.
\begin{enumerate}
\item There exist finite $K(A,1)$-spaces. \cite{CD2}, \cite{BrWa}
\item $A$ is torsion-free.  \cite{BS} 
\item The center of $A$ is infinite cyclic. \cite{BS}
\item $A$ is linear.  \cite{Big}, \cite{Kra}, \cite{CW}
\item  $A$ is biautomatic and hence has solvable word and conjugacy problem.   \cite{Del}, \cite {BS}, \cite{Ch1}, \cite{Ch2}
\end{enumerate}

By contrast, none of the above properties are known to hold for general, non-spherical Artin groups.  For some partial results on non-spherical Artin groups see \cite{Alt},  \cite{BMc}, \cite{Ch3}, \cite{CD1}, \cite{CD2}, \cite{ChP}, \cite{Che}, \cite{Dig}, \cite{Pei}.

One question which is still open for both spherical and non-spherical Artin groups is whether they are CAT(0) groups, that is, whether they act properly, cocompactly by isometries on a CAT(0) space (see below for definition).  For partial results on this problem see  \cite{Bel}, \cite{Bes}, \cite{BrCr} and \cite{TBr}.


\section{Right-angled groups and CAT(0) cube complexes}

\subsection{Right-angled Coxeter and Artin groups}

For the remainder of this paper, we will concentrate on a very special class of non-spherical Artin groups known as right-angled Artin groups.  A \emph{right-angled Artin group} is one in which 
$m_{ij} \in \{2,\infty\}$ for all $i,j$.  In other words, in the presentation for the Artin group, all relations are commutator relations:  $s_is_j=s_js_i$.  Right-angled Coxeter groups are defined similarly. 

The easiest way to specify the presentation for a right-angled Coxeter or Artin group is by means of the \emph{defining graph} $\G$.  This is the graph whose vertices are labeled by the generators $S=\{s_1, \dots ,s_n\}$ and whose edges connect a pair of vertices $s_i,s_j$ if and only if $m_{ij}=2$.  Note that \emph{any} finite, simplicial graph $\G$ is the defining graph for a right-angled Coxeter group $W_\G$ and a right-angled Artin group $\AG$.

\begin{figure}
\begin{center}
\includegraphics[height=2cm]{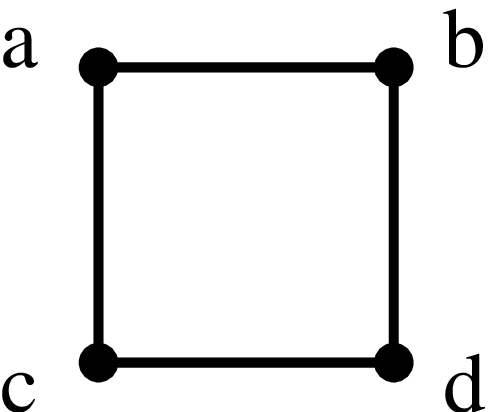}
\end{center}
\caption{$\AG=F(a,c)\times F(b,d)$ } 
\label{Fig2}
\end{figure}

\begin{examples} At one extreme, we have a graph  $\G$ with no edges, in which case $\AG=F_n$, the free group on $n$ generators. At the other extreme is the case of a complete graph $\G$  (every vertex is connected to every other vertex by an edge), in which case $\AG=\Z^n$, the free abelian group on $n$ generators.  Between these are many possibilities.  For example, if $\G$ is square as in Figure \ref{Fig2},  then $\AG$ decomposes as a direct product of two free groups $\AG=F(a,c) \times F(b,d)$.  Whereas if $\G$ is an $n$-gon for $n \geq 5$, then $\AG$ cannot be decomposed as either a direct product or a free product. Note that the only case in which $\AG$ is a spherical Artin group is when $\G$ is a complete graph.
\end{examples}

While any Artin group has a Coxeter group quotient, right-angled Artin groups also embed as subgroups of Coxeter groups. Indeed, Davis and Januszkiewicz \cite{DJ} prove:

\begin{theorem}
Every right-angled Artin group embeds as a finite index subgroup of a right-angled Coxeter group.
\end{theorem}

It follows, in particular, that right-angled Artin groups are linear.  (Linearity was also proved by Hsu and Wise in \cite{HW1}.)  For an idea of the proof of this theorem, consider the Coxeter group $D_\infty=\Z/2\Z \ast \Z/2\Z$ whose defining graph consists of two disjoint points.  $D_\infty$ acts on the real line as the group generated by (affine) reflections across two distinct points. The one generator Artin group $\Z$ embeds in $D_\infty$  as the index 2 subgroup of translations.  For a general right-angled Artin group $\AG$,  Davis and Januszkiewicz define a ``doubling"  $\tilde\G$ of $\G$  such that each vertex  $s_i$ of $\G$ corresponds to a pair of vertices $s_i, s'_i$ in $\tilde\G$. The new vertices $s_i'$ are connected to every other vertex except $s_i$ by an edge. The inclusion $\AG \to W_{\tilde\G}$ is defined by sending $s_i$ to the ``translation" $s_is'_i$.

Right-angled Artin groups also arise as subgroups of arbitrary Artin groups.  It was conjectured by Tits and proved by Crisp and Paris \cite{CrPa} that for any Artin group $A$, the subgroup generated by the squares of the generators, $s_1^2, \dots ,s_n^2$ is a right-angled Artin group; the only relators are commutators between squares of generators that commute in the original group $A$.

\subsection{Subgroups}

For any subset $T \subset S$, we denote by $A_T$  the subgroup of $\AG$  generated by $T$.  This subgroup is isomorphic to the Artin group associated to the subgraph of $\G$ spanned by $T$.  Likewise, $W_T$ denotes the subgroup of $W$ spanned by $T$. A subgroup of the form $A_T$ (or $W_T$) is called a \emph{special subgroup}. 

Right-angled Artin groups contain many interesting subgroups.   Droms \cite{Dro1} shows $\AG$ is a 3-manifold group if and only if every component of $\G$ is a tree or a triangle. In the tree case, $\AG$ is the fundamental group of a graph manifold, that is, a 3-manifold which decomposes along embedded tori and Klein bottles into finitely many Siefert fiber spaces.  More general right-angled Artin groups contain such 3-manifold groups as special subgroups.   Crisp and Wiest \cite{CrWi}, extending work of Droms, Servatius, and Servatius \cite{DSS}, show that with only 3 exceptions, any surface group embeds in a right-angled Artin group.  They also show that ``graph braid groups" embed in right-angeld Artin groups. These are discussed in Section \ref{robots} below.  More exotic subgroups of right-angled Artin groups, introduced by Bestvina and Brady \cite{BB} are discussed in Section \ref{Bestvina-Brady}.

The question of subgroup separability for right-angled Artin groups has been studied by Hsu and Wise  \cite{HW2} and by Metafsis and Raptis  \cite{MeRa1}, \cite{MeRa2}.  A subgroup $H$ of a group $G$ is \emph{separable} if for every element $g \in G \backslash H$, there is a finite quotient of $G$ to which $g$ maps non-trivially.  The group $G$ is \emph{subgroup separable} if every finitely generated subgroup of $G$ is separable.  Metafsis and Raptis prove,

\begin{theorem}
$\AG$ is subgroup separable if and only if $\G$ has diameter at most 2 and contains no square as a full subgraph.
\end{theorem} 

\subsection{The word problem}  

Given a presentation $\langle S ~|~ R \rangle$  for a group $G$, there are, in general, many ways of representing an element $g$ as a word in the generators $S$.  The \emph{word problem} for $G$ asks whether there is an algorithm for deciding when two words in $S$ represent the same element of $G$.  Solutions to the word problem often involve finding a canonical form (i.e., a preferred word) to represent each element of $G$.  A particularly nice situation is when $G$ has a \emph{biautomatic structure} in which the word problem is solvable by means of a collection of finite state automata (see \cite{WP} for a precise definition).  In this case, there is also a solution to the conjugacy problem, that is, there is an algorithm to determine when two words represent conjugate elements of $G$.

One of the key differences between spherical and non-spherical Artin groups is that for the former, there exist nice canonical forms for elements of $A$ as words in the standard generating set $S$.  This canonical form was introduced by Garside \cite{Gar} for the braid groups and extended to spherical Artin groups by Deligne \cite{Del} and Brieskorn-Saito \cite{BS}. They give rise to a biautomatic structure on $A$ \cite{Ch1},\cite{Ch2} and allow for easy solutions to many algebraic questions.  By contrast, for most non-spherical Artin groups, no algorithmic solution to the word problem is known, but in the case of right-angled Artin groups,  there is such a solution which we now describe.  

Consider the set $Min(a)$ of all minimal length words in $S$ representing a fixed element $a \in \AG$.  In her thesis \cite{Gre}, Green proves that the set $T \subset S$ of all elements which occur as the initial letter in such a word is a spherical subset of $S$, that is, the elements of $T$ mutually commute.  (This can be also be proved geometrically using the action of $\AG$ on a CAT(0) cube complex described in the next section.)  It follows that there is a unique maximal length element $w_0$ in the abelian group $A_T$ which can occur as an initial segment in $Min(a)$.  Writing $a$ as a product $a=w_0a_1$, we can now repeat this process starting with $a_1$ to get a maximal initial segment $w_1$ in $Min(a_1)$ and a decomposition $a=w_0w_1a_2$.  Continuing this process, we eventually end up with a decomposition of $a$ into $a=w_0w_1 \dots w_k$ with each $w_i$ lying in an abelian special subgroup $A_{T_i}$.  We call this the canonical form for $a$.  It is unique up to the commutation relations in $A_{T_i}$.  (The reader familiar with Garside normal forms for spherical  Artin groups will note the similarity.)

In \cite{HM}, Hermiller and Meier prove that, starting with any word representing $a$, the normal form for $a$ can be obtained by a process of  shuffling commuting elements and canceling inverse pairs of generators whenever possible.  In particular, the normal form can be obtained by applying operations to a word that do not increase its length. (More formally, they obtain a ``rewriting system" for $\AG$.)  They also prove that this normal form gives a biautomatic structure on $\AG$.  An alternate proof of biautomaticity is given by VanWyk in \cite{vW}.   

Similar arguments can be used to determine centralizers of elements in $\AG$.  In \cite{Ser}, Servatius proves that any conjugacy class of elements in $\AG$ contains an essentially unique (up to cyclic permutation) element of minimal length. 
If $u$ is such a minimal length element, and $A_T$ is the smallest special subgroup containing $u$, let $A_T=A_{T_1} \times \dots \times A_{T_k}$ be the maximal decomposition of $A_T$ as a direct product special subgroups (i.e., the factors cannot be further decomposed).  Then $u$ can be written as a product $u=u_1^{m_1} \dots u_k^{m_k}$ where $u_i \in A_{T_i}$ and $m_i$ is maximal.  Servatius proves that the centralizer of $u$ is
\[
C(u)=\langle u_1\rangle \times \dots \times \langle u_k\rangle \times \langle link(T)\rangle
\]
where $link(T)$ is the set of vertices $v \notin T$ which are adjacent to every vertex in $T$.
In particular, for a single vertex $s$, $C(s)= \langle s\rangle  \times \langle link(s) \rangle$, the special subgroup generated by $s$ and the vertices adjacent to $s$. Hence, the center of $\AG$ is the special subgroup generated by the (possibly empty) set of vertices in $\G$ that are connected to every other vertex by an edge.

\subsection{CAT(0) cube complexes}

Associated to right-angled Coxeter and Artin groups are certain CAT(0) cube complexes which play a central role in applications.  We review very briefly the basics of CAT(0) cube complexes and refer the reader to Bridson and Haefliger's book \cite{BH} for more details. 

A metric space $X$ is \emph{geodesic} if any two points are connected by a geodesic segment, i.e., an isometric embedding of an interval into $X$. It follows that the distance between any two points is equal to the length of the shortest path between them.  A geodesic metric space is \emph{CAT(0)} if geodesic triangles in $X$ are at least as thin as comparison triangles of the same side lengths  in the Euclidean plane.  That is, the distance between any two points on the triangle in $X$ is less than or equal to the Euclidean distance between the corresponding points on the comparison triangle. A space which is locally CAT(0) is said to have \emph{non-positive curvature}.  

A fundamental theorem in CAT(0) geometry states that a space which is locally CAT(0) (non-positively curved) and simply connected, is globally CAT(0).  It is easy to see that a  CAT(0) space is contractible, hence any geodesic metric space of non-positive curvature is aspherical. 
 
One way to construct spaces of non-positive curvature is to begin with a cubical complex $K$ (the cubical analogue of a simplicial complex).  Assign each cube the metric of a regular Euclidean cube of side length $1$, and take the induced geodesic metric on $K$.  We will refer to this as the \emph{cubical metric} on $K$.  While in general it is very difficult to tell whether a geodesic metric has non-positive curvature, for cubical metrics, one has a simple, combinatorial condition on links of vertices which determines whether the complex is locally CAT(0).  The \emph{link} of a vertex $v$ in $K$ is the simplicial complex with one $k$-simplex for each $(k+1)$-cube containing $v$. Geometrically, one can think of the link of $v$ as the boundary of a small ball centered at $v$.  

A simplicial complex $L$ is called \emph{flag} if every complete subgraph in the 1-skeleton spans a simplex in $L$.  In other words, $L$ is completely determined by its 1-skeleton.  In \cite{Gro}, Gromov proves the following theorem.

\begin{theorem}\label{flag}
The cubical metric on $K$ has non-positive curvature if and only if the link of every vertex is a flag complex.
\end{theorem}

CAT(0) cubical complexes also have the advantage that a group $G$ acting properly, cocompactly by isometries on such a complex has more structure than a general CAT(0) group. For example, $G$ is biautomatic \cite{NR} and satisfies a Tit's alternative: any subgroup of $G$ is either virtually abelian or contains a non-abelien free group \cite{SaWi}.

\subsection{The Davis complex}

Associated to any Coxeter group $W$ is a contractible, polyhedral cell complex $\Sigma_W$ called the \emph{Davis complex}, introduced by M.~Davis in \cite{Dav}.  In the case of a right-angled Coxeter group $W_\G$,  we will denote the Davis complex  by $\Sigma_\G$.  It is the cubical complex whose 1-skeleton is the Cayley graph of $W_\G$ and a set of vertices $C \subset W_\G$  spans a cube if and only if it forms a coset $C=wW_T$ for some \emph{finite} special subgroup $W_T$.  (Here, we allow $T=\emptyset$, in which case $wW_T=\{w\}$ is a vertex.)
$W_\G$ acts on this complex (by left multiplication) with finite stabilizers.  Thus, any torsion-free subgroup of finite index acts freely on $\Sigma_\G$.

\begin{figure}
\begin{center}
\includegraphics[height=4cm]{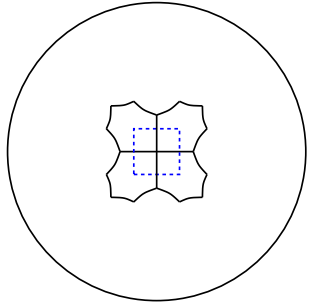}
\end{center}
\caption{Tiling of $\Hy^2$ by right-angled pentagons} 
\label{Fig3}
\end{figure}

A good example to consider is where the defining graph $\G$ is a pentagon.  In this case, $W_\G$ can be realized as the hyperbolic reflection group generated by reflections in the sides of a right-angled pentagon $P$ in $\Hy^2$.  Translates of $P$ by the group give a tiling of $\Hy^2$ with four pentagons meeting at each vertex and $\Sigma_\G$ can be identified (as a cell complex) with the dual cellulation of $\Hy^2$ by squares. (Figure \ref{Fig3})

Since $W$ acts transitively on the vertices of $\Sigma_\G$,  the link of every vertex is isomorphic.  This link, which we denote by $N_\G$, is called the \emph{nerve of $W$}. It plays an important role in many applications, so let's take a closer look at it. For simplicity, we choose our vertex to be $W_\emptyset $.  Then a $k$-simplex in $N_\G$ corresponds to a $(k+1)$-cube of $\Sigma_\G$  containing the vertex $W_\emptyset$.  Such a cube is spanned by the vertices in a finite special subgroup $W_T$ with $|T|=k+1$.   A set $T\subset S$ generates a finite subgroup of $W_\G$  if and only if  $T$ spans a complete subgraph in $\G$.  It follows that  $N_\G$ is the simplicial complex consisting of a simplex  for each complete subgraph  in $\G$.   In other words, the nerve $N_\G$ is the flag complex whose 1-skeleton is $\G$.

In particular, the link of every simplex in the Davis complex is a flag complex.  It is not difficult to show that the Davis complex is simply connected.  Applying Theorem \ref{flag}, we thus obtain

\begin{theorem}
For any $\G$, the cubical metric on $\Sigma_\G$ is CAT(0).  
\end{theorem}

\subsection{The Salvetti complex}

Associated to a right-angled Artin group $\AG$ is an analogous cube complex $\tilde\Sa$ constructed as follows.  Begin with a wedge of circles attached to a point $x_0$ and labeled by the generators $s_1, \dots ,s_n$.  For each edge, say from $s_i$ to $s_j$ in $\G$, attach a 2-torus with boundary labeled by the relator $s_is_js_i^{-1}s_j^{-1}$.  For each triangle in $\G$ connecting three vertices $s_i,s_j,s_k$, attach a 3-torus with faces corresponding to the tori for the three edges of the triangle.  Continue this process, attaching a $k$-torus for each set of $k$ mutually commuting generators (i.e., generators spanning a complete subgraph in $\G$).  The resulting space, $\Sa$,  is called the \emph{Salvetti complex} for $\AG$. It is clear, by construction, that the fundamental group of $\Sa$ is $\AG$.\begin{footnote}{The term ``Salvetti complex" is more often used to denote the cover of $\Sa$  whose fundamental group  is the kernel of the projection $\AG \to W_\G$. Since we are not concerned with this cover, we will abuse the terminology and call $\Sa$ the Salvetti complex.  The analogue of this complex for spherical Artin groups was introduced by Salvetti in \cite{Sal}. }\end{footnote}  

\begin{figure}
\begin{center}
\includegraphics[height=2.2 cm]{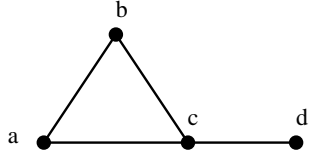}
\end{center}
\caption{A defining graph $\G$} 
\label{Fig4}
\end{figure}

For example, the Salvetti complex for a pentagon is made up of five 2-tori, each glued to the next along one curve, while the Salvetti complex for the graph $\G$ in Figure \ref{Fig4} is made up of one 3-torus and one 2-torus glued along a single curve corresponding to the vertex $c$.

There is also a description of $\tilde\Sa$ which is analogous to that of $\Sigma_\G$.  Namely, for any finite special subgroup $W_T$ of the Coxeter group, with $T=\{t_1, \dots ,t_k\}$, we can lift $W_T$ to the subset (not subgroup) $\widehat W_T \subset A_T$ consisting of elements of the form $t_1^{\epsilon_1}\cdots  t_k^{\epsilon_k}$ with $\epsilon_i=0,1$.  Then $\tilde\Sa$ is the cube complex whose vertices are the elements of $\AG$ and whose cubes are sets of vertices of the form $a\widehat W_T$ for some $a \in \AG$ and some finite special subgroup $W_T$.

The following theorem is proved in \cite{CD2}.

\begin{theorem} The universal cover of the Salvetti complex,
$\tilde \Sa$, is a CAT(0) cube complex. Hence $\Sa$ is a $K(\AG,1)$ space.
\end{theorem}

It is clear by construction that the fundamental group of $\Sa$ is $\AG$, thus to prove the theorem we need only check that links in $\tilde\Sa$ are flag complexes.  All links of vertices in $\tilde\Sa$ are equivalent to the link $L_\G$ of the unique vertex $x_0$ in $\Sa$.  This link consists of a $(k-1)$-sphere for each $k$-torus in $\Sa$.  As a simplicial complex, it has two vertices $s_i^+$ and $s_i^-$ for each vertex $s_i$ in $\G$ (corresponding to the incoming and outgoing tangent vectors to the loop labeled $s_i$).  A set of vertices $s^{\pm}_{i_1}, \dots s^{\pm}_{i_k}$ spans a simplex in $L_\G$ if and only if the $s_{i_j}$'s are distinct and span a complete subgraph of $\G$, or in other words, if and only if they represent pairwise commuting generators in $\AG$. (Note that in this case, the set of \emph{all} simplices  $s^{\pm}_{i_1}, \dots s^{\pm}_{i_k}$  forms a $(k-1)$-sphere.)
From this description, it is easy to see that $L_\G$ is a flag complex.

It follows from the theorem that $\AG$ has a finite $K(\pi,1)$-space and hence is torsion-free.  It also follows that $\AG$  acts freely, cocompactly on a finite dimensional CAT(0) space and so, by the work of Niblo and Reeves \cite{NR},  we obtain another proof that $\AG$ is biautomatic.   
Other algebraic properties of right-angled Artin groups arising from their actions on CAT(0) cube complexes have been explored by Hsu and Wise in \cite{HW2}. 

\subsection{Cohomology and the topology at infinity}

One use for a $K(G,1)$-space $X$ is to determine the cohomology of the group $G$.  In the case of a right-angled Artin group, the cohomology can be easily computed from its Salvetti complex.  Indeed, it is easy to see that every cell in the Salvetti complex is a cycle (the cells are glued on as tori), so $H_k(\AG)$ is a free abelian group with one generator for each $k$-cell in $\Sa$, or equivalently, one generator for each $(k-1)$-cell in the nerve $N_\G$.  

We can also use the Salvetti complex to determine the ``topology at infinity"  of $\AG$.  The topology at infinity of a group $G$ with a finite $K(G,1)$-space $X$ is determined by considering complements of larger and larger compact sets in the universal covering space $\tilde X$.   It gives rise to useful asymptotic invariants of the group.  For example, the number of ``ends" of a group is the supremum over compact sets $C$  of the number of connected components in $\tilde X \backslash C$.  Thus, a group which is ``connected at infinity" is a 1-ended group and hence cannot be decomposed non-trivially as an amalgamated product or HNN extension over a finite group.  In \cite{BrMe}, Brady and Meier prove that $\AG$ is $m$-connected (respectively $m$-acyclic) at infinity if and only if  $L_\G$ is $m$-connected (respectively $m$-acyclic).  (This result can also be reformulated in terms of links in $N_\G$.)

Topology at infinity can also be determined from the cohomology of  $G$ with coefficients in the group ring $\Z G$.  For example, a group is 1-ended if and only if $H^1(G, \Z G)$ is trivial.  The cohomology of $G$ with group ring coefficients is the same as the cohomology with compact supports of $\tilde X$.   For Coxeter groups, the cohomology with group ring coefficients was computed by Davis in \cite{Dav2}.  In \cite{JeMe}, Meier and Jensen use the fact that right-angled Artin groups embed as finite index subgroups of right-angled Coxeter groups to compute $H^*(\AG;\Z\AG )$ for any right-angled Artin group. The answer is given in terms of the (integral) cohomology of links of simplices in $N_\G$.  One corollary of their computation is a characterization of when $\AG$ is a duality group or Poincar\'e duality group.  A group $G$ is said to be a \emph{duality group} if there is a module $D$ and a positive integer $n$ such that $H_i(G;M) \cong H^{n-i}(G; M \otimes D)$ for all $i$ and all $G$-modules $M$.  It is a \emph{Poincar\'e duality group} if, in addition, $H^n(\AG;\Z\AG ) =\Z$.  Meier and Jensen's computation shows that $\AG$ is a duality group if and only if $N_\G$ is Cohen-Macauley (i.e. the cohomology of the link of every simplex is torsion-free and concentrated in the top dimension), and $\AG$ is a Poincar\'e duality group if and only if $\AG$ is free abelian.  A different proof of this characterization appears in \cite{BrMe}.

Another topic of current interest is the $\ell_2$-cohomology of Coxeter and Artin groups.  The $\ell_2$-cohomology of $G$ is computed from the complex of square-summable cochains on $\tilde X$.  The $\ell_2$ cohomology of right-angled Coxeter groups is studied by Davis and Okun in \cite{DO} and of right-angled Artin groups by Davis and Leary in \cite{DL}.  Davis and Leary prove (as a special case of a more general theorem) that the reduced $\ell_2$-cohomology of $\AG$ is isomorphic to $\bar H^*(N_\G) \otimes \ell_2(\AG)$ where $\bar H^*$ denotes ordinary (reduced) cohomology and $\ell_2(\AG)$ is the Hilbert space of square-summable real-valued functions on $\AG$.

Other cohomological invariants of $\AG$, such as lower central series quotients, are computed by Papadima and Suciu in \cite{PaSu}.


\section{Automorphism and quasi-isometries of right-angled Artin groups}

Recently, attention has turned to understanding automorphisms and quasi-isometries of right-angled Artin groups. 
Recall that a \emph{quasi-isometry} between two metric spaces $X$ and $Y$ is a map $f: X \to Y$ for which there exist constants $c \geq 0, k \geq 1$ such that
\[
\frac{1}{k}\,d(f(x_1),f(x_2)) -c  \leq d(x_1,x_2) \leq k \,d(f(x_1),f(x_2)) +c 
\]
 for all $x_1,x_2 \in X$, and every  point in $y$ lies in the $c$-neighborhood of $f(X)$. Quasi-isometry is particularly well-suited to studying geometric properties of finitely generated groups since the word metrics with respect to different generating sets are quasi-isometric.

In \cite{Dro2},  Droms proved that two right-angled Artin groups $A_\G$ and $A_{\Gamma'}$ are isomorphic as groups if and only if their defining graphs are isomorphic. Viewing the groups as metric spaces (via the word metric), one can ask for a weaker classification, namely which right-angled Artin groups are \emph{quasi-isometric} to each other?   Some partial answers to these questions have recently been found by Behrstock and Neumann \cite{BeNe} and by Bestvina, Kleiner, and Sageev \cite{BKS}.   Their results suggest a fairly complicated picture.  On the one hand, Behrstock-Neumann show that many distinct graphs give rise to quasi-isometric Artin groups (right-angled and otherwise), while Bestvina-Kleiner-Sageev show that other  right-angled Artin groups are very rigid.  For example,  \emph{all} right-angled Artin groups whose defining graphs are trees of diameter at least 3 (as well as some non-right-angled Artin groups) belong to a single quasi-isometry class.  
By contrast,
for groups whose defining graph is ``atomic"  (no valence 1 vertices, no cycles of length $<5$ and no separating vertex stars),  $\AG$ is quasi-isometric to $A_{\G'}$ if and only if $\G \cong\G'$, and moreover, any isomorphism from $\AG$ to $A_{\G'}$ is induced by an isometry from $\tilde\Sa$ to $\tilde S_{\G'}$.

Another interesting area to explore is the automorphism group of a right-angled Artin group. 
For a group $G$, we denote by $Aut(G)$ the automorphism group of $G$ and by $Out(G)$ the outer automorphism group, that is, the quotient of $Aut(G)$ by the group of conjugation automorphisms.  Automorphism groups play an important role in algebra and topology.  For example, the mapping class group of a closed surface $S_g$ of genus $g$ is the outer automorphism group of  $\pi_1(S_g)$.  Two automorphism groups which have been particularly important are $GL(n,\Z)$, the automorphism group  of a free abelian group, and $Aut(F_n)$, the automorphism group of a free group.  Although Formanek and Procesi \cite{FP} have shown that $Aut(F_n)$ is not a linear group for $n\geq 3$, these two classes of groups have many properties in common.  Since right-angled Artin groups ``interpolate" between free groups and free abelian groups, it is natural to ask which properties hold for automorphism groups of all right-angled Artin groups. 

Little is known about  the automorphism groups of  right-angled Artin groups in general with the exception of the work of Servatius \cite{Ser} and Laurence \cite{Lau} who describe a finite generating set for these groups.   In joint work of the author, John Crisp and Karen Vogtmann \cite{CCV}, we study the automorphism group of 2-dimensional right-angled Artin groups ($\AG$ is \emph{2-dimensional} if its Salvetti complex $\Sa$ is 2-dimensional, or equivalently, if $\G$ has no triangles).  Assume for the remainder of this section that $\G$ is connected and triangle-free.

A set of vertices $J$ in $\G$ is called a \emph{join} if it spans a complete bipartite subgraph.  That is, $J$ is the disjoint union of two non-empty subsets $V_1$ and $V_2$ such that every vertex in $V_1$ is connected to every vertex in $V_2$.   The no-triangle condition then guarantees that no two vertices in the same $V_i$ are connected by an edge.  It follows that the special subgroup $A_J$ generated by $J$ is a direct product of free groups, $A_J=F(V_1) \times F(V_2)$.   A join $J$ is a \emph{maximal join}  if it is not properly contained in another join. 

The maximal joins provide the key structure for $Out(\AG)$.  Up to diagram symmetry, automorphisms preserve the maximal join subgroups $A_J$ up to conjugacy.  This gives rise to a homomorphism
\[
\phi: Out^0(A_{\Gamma}) \to \prod Out(A_J)
\]
where $Out^0(A_{\Gamma})$ is a finite index normal subgroup of $Out(A_{\Gamma})$ which avoids certain diagram symmetries and the product is taken over maximal joins in $\G$.  We prove,  

\begin{theorem}\label{kernel}
For a vertex $s$ in $\G$, let $k_s$ denote the number of connected components of $\G -\{ s\}$ containing at least 2 vertices.
Then the kernel $K$ of $\phi$ is a finitely generated, free abelian group of rank $\sum_{s \in S} (k_s -1)$.
\end{theorem}

It follows from this theorem that $Out(\AG)$ contains a torsion-free subgroup of finite index.  Moreover, $K$ commutes with the group of "leaf transvections" (multiplication of a leaf vertex by the adjacent vertex)  and together they generate a free ableian subgroup of $Out (\AG)$ of rank 
$\sum_{s \in S} (m_s -1)$  where $m_s$ denotes the number of \emph{all} connected componenets of $\G -\{ s\}$.

A key tool in the study of  automorphism groups of free groups is their action on a topological space  $\OF$ known as \emph{outer space}.  Introduced by Culler and Vogtmann in \cite{CV}, outer space is a finite dimensional, contractible space with a proper action of $Out(F_n)$.  Points in this space correspond to (free, minimal, isometric) actions of $F_n$ on a tree.  To obtain an analogous space $\mathcal O(A_{\Gamma})$ for a 2-dimensional Artin group $\AG$, 
it is natural to consider actions of $\AG$ on CAT(0) 2-complexes. 

Consider, for example, the 
case in which the defining graph consists of a single join $J$. In this case, $A_J$ is a product of two free groups and any (free, minimal, isometric) action on a CAT(0) 2-complex $X_J$ gives rise to a splitting of $X_J$ as a product of two trees $X_J=T_1 \times T_2$.   (It is important to note, however, that if one of the free factors is cyclic, then one of the trees is just a line and the action may have shearing; the two free factors need not act orthogonally.)

For a more general defining graph, an action of $\AG$ on a CAT(0) 2-complex $X$ gives rise to a covering of $X$ by such tree products $X_J$, where $J$ ranges over maximal joins in $\G$.  Composing the action on $X$ with an automorphism of $\AG$ gives rise to a new action, and hence a new set of subspaces $X'_J$. These subspaces (together with their $A_J$-actions) are sufficient to capture the action $Out(\AG)$ on $\AG$.  We use these (tree $\times$ tree)-spaces to construct an analogue of outer space for $\AG$ and we prove,  

\begin{theorem}
For any connected, triangle-free graph  $\Gamma$,  outer space  $\mathcal O(A_{\Gamma})$ is finite dimensional, contractible, and has a proper action of $Out^0(A_{\Gamma})$. 
\end{theorem}

The dimension of $\mathcal O(A_{\Gamma})$ is large, but as in the case of a free group, it retracts onto a ``spine" of lower dimension. The dimension of the spine provides an upper bound on the virtual cohomological dimension and the rank of an abelian subgroup provides a lower bound. 

\begin{theorem} Let $\G$ be connected and triangle-free, let $\ell_s$ be the valence of a vertex $s$ in $\G$ and let $m_s$ be the number of connected components in $\G-\{s\}$. Then $Out(\AG)$ is virtually torsion-free and its virtual cohomological dimension satisfies
 \[
\sum_{s \in S} (m_s -1)  \leq vcd (Out(\AG)) <  \sum_{s\in S}3( \ell_s-1).
\]
\end{theorem}

In the case that $\G$ is a tree with $e$ edges, $m_s=l_s$ and $\sum_{s\in S}( \ell_s-1)=e-1$ so the vcd satisfies $e-1 \leq vcd (Out(\AG)) < 3(e-1)$.



\section{Applications of right-angled Coxeter and Artin groups}\label{applications}

Right-angled Coxeter and Artin groups appear in many applications.  These groups, and the CAT(0) complexes associated to them, have provided counter-examples to a number of conjectures in topology and group theory.  The key to these counter-examples is the ability to design the nerve $N_\G$ to serve a particular purpose.  Note that the barycentric subdivision, $K'$ of any simplicial complex $K$ is a flag complex.  Taking $\G$ to be the 1-skeleton of $K'$, the associated nerve, $N_\G$, is precisely $K'$.  In particular, $N_\G$ can be chosen to be homeomorphic to any finite simplicial complex.  This was first used by Davis \cite{Dav} in his well-known construction of closed, aspherical manifolds not covered by Euclidean space.

\subsection{The Davis manifolds}  

In \cite{Dav}, Davis introduced the complex $\Sigma_\G$ and observed that one could construct many aspherical manifolds using this complex. 
A neighborhood of a vertex in $\Sigma_\G$ is homeomorphic to the cone on $N_\G$ and these cones are fundamental domains for the $W_\G$-action (they are dual to the cubical structure).  In particular,  if we take $N_\G$ to be an $(n-1)$-sphere, then  $\Sigma_\G$ is an $n$-manifold and this manifold is contractible since $\Sigma_\G$ supports a CAT(0) metric.  $W_\G$ acts cocompactly as a discrete reflection group on this manifold, hence for any finite-index, torsion-free subgroup $G \subset W_\G$,  the orbit space $\Sigma_\G /G$ is a closed, aspherical manifold.  

In \cite{Dav}, Davis uses a variation on this construction to produce  examples of closed, aspherical manifolds of dimension $\geq4$ that are not covered by Euclidean space.  His construction is roughly as follows. If $N_\G$ is a generalized homology sphere instead of an actual sphere,  Davis shows that one can replace a  neighborhood of each vertex in $\Sigma_\G$ with a compact manifold with boundary to obtain a contractible manifold $M$  with a cocompact $W_\G$-action.  He proves that the resulting manifold $M$ is simply connected at infinity if and only if the nerve, $N_\G$, was simply connected.  Taking $N_\G$ to be a non-simply-connected generalized homology sphere,  and $G \subset W_\G$ to be a torsion-free subgroup of finite index, gives the desired examples $M/G$. 

The idea of constructing manifolds with specified links was also used by the author and Davis in their work on the combinatorial Hopf conjecture.  An old conjecture of Hopf states that the Euler characteristic of a $2n$-dimensional closed, Riemannian manifold $M$ of non-positive sectional curvature should satisfy $(-1)^n\chi (M)\geq 0$.  In \cite{CD3}, the author and M.~Davis considered the analogous conjecture for manifolds with a cubical metric of non-positive curvature, i.e., cube complexes whose links are flag triangulations of a $(2n-1)$-sphere.  In this case, the Euler characteristic of $M$ can be written as a sum of ``local Euler characteristics" at each vertex, in which the contribution of a $k$-cube is divided among its $2^{k}$ vertices.  The local Euler characteristic at $v$ depends only on the link $L_v=link(v,M)$:
\[
\chi(v,M)=1+\sum_{\sigma \in L_v} \left(\frac{-1}{2}\right)^{\dim(\sigma)+1}
\]
We conjectured that each of these \emph{local} Euler characteristics must  satisfy, $(-1)^n\chi(v,M)  \geq 0$ (and showed that this was the case under some additional hypotheses).  Clearly this would imply the Hopf conjecture for cubical manifolds.  It is less apparent that the converse is true. However, given a flag triangulation $L$ of the $(2n-1)$-sphere, we can take $\G$ to be the one skeleton of $L$, and $G$ to be a torsion-free subgroup of finite index in $W_\G$.  This gives a closed, cubical manifold of non-positive curvature,  $M=\Sigma_\G /G$, all of whose links are precisely $L$.  It follows that any counterexample to this local Hopf conjecture gives a counterexample to the global Hopf conjecture for cubical manifolds.  For progress on these and related conjectures, see \cite{DO}, \cite{Bra}, and \cite{ReWe}.

\subsection{The Bestvina-Brady groups}\label{Bestvina-Brady}

The flexibility of $N_\G$ is also central to Bestvina and Brady's remarkable results on 
finiteness properties of groups.  In working with infinite groups, we can ask that the group satisfy various finiteness conditions.  The most well known of these are the conditions that the group $G$ be finitely generated or finitely presented.  These conditions can be reinterpreted as saying that there exists a $K(G,1)$-space with a finite 1-skeleton, respectively, a finite 2-skeleton.  Thus, a natural generalization is to ask that there exist a $K(G,1)$ space with a finite $n$-skeleton. In this case we say that $G$ is of type $F_n$.  Finally, we might ask that $G$ have a $K(G,1)$-space which is a finite complex, in which case we say that $G$ is of type $F$.  

There are also homological versions of these conditions.  Note that if $X$ is a $K(G,1)$-space, then the cellular chain complex of its universal cover, $C_\ast(\tilde X)$ is an exact sequence of projective
$\Z G$-modules.  Thus, if $G$ is of type $F_n$, then there exists an exact sequence (or ``resolution") of finitely generated projective $\Z G$-modules,
\[
P_n \to  P_{n-1} \to \dots \to P_0 \to \Z \to 0.
\]
A group with such a resolution is said to have type $FP_n$.  If this resolution terminates, that is, if there exists a finite resolution 
\[
0\to P_n \to  P_{n-1} \to \dots \to P_0 \to \Z \to 0,
\]
for some $n$, then $G$ is said to have type $FP$.  In particular, any group of type $F$ is also of type $FP$.

Until recently, it was unknown whether the homological conditions $FP_n$ (respectively $FP$) implied the geometric conditions $F_n$ (respectively $F$) for $n>1$.  In \cite{BB}, M.~Bestvina and N.~Brady show that this is not the case.  Their counterexamples are constructed using right-angled Artin groups and Salvetti complexes. For any right-angled Artin group we can define a homomorphism  $\phi: A_\G \to \Z$ by sending each generator $s_i$ to $1 \in \Z$.  This map can be realized as a map of $K(\pi,1)$-spaces, $S_\G \to \mathbb S^1$.  Lifting to the universal cover gives an equivariant ``Morse function" $f: \tilde S_\G \to \R$.  

The kernel $K_\G$ of $\phi$ acts on the level sets $f^{-1}(r)$ of this Morse function.  By analyzing the topology of the level sets, Bestvina and Brady are able to determine finiteness properties of the group $G$.  Since $\tilde S_\G$ is contractible, the topology of the level set can be determined from the topology of the half spaces $f^{-1}(-\infty, r]$ and $f^{-1}[r, \infty)$ via the Mayer-Vietoris sequence.  These in turn, are determined by the ``upward" and ``downward" links of vertices in $\tilde S_\G$.  The key to their construction, is that these upward and downward links turn out to be precisely the complexes $N_\G$ which arose in the Davis construction.  In particular, they can be chosen to be homeomorphic to any finite simplicial complex.  Bestvina and Brady prove,

\begin{theorem} Let $K$ be the kernel of the homomorphism $\phi: A_\G \to \Z$. Then
\begin{enumerate}
\item $K$ is type $FP_n$ if and only if $H_i(N_\G)=0$ for all $i<n$,
\item $K$ is type $FP$ if and only if  $N_\G$ is acyclic,
\item $K$ is finitely presented (type $F_2$) if and only if $N_\G$ is simply connected.
\end{enumerate}
\end{theorem}

In particular, choosing $\G$ such that $N_\G$ is acyclic but not simply connected gives an example of a group which is type $FP$ (and hence $FP_n$ for all $n$) but not type $F_2$ (and hence not $F_n$ for any $n \geq 2$).  

For more on the Bestvina-Brady groups see  \cite{MMW}, \cite{MvW}, \cite{LeNu}, and \cite{PaSu}.

\subsection{The Croke-Kleiner examples}

\begin{figure}
\begin{center}
\includegraphics[height=2.2 cm]{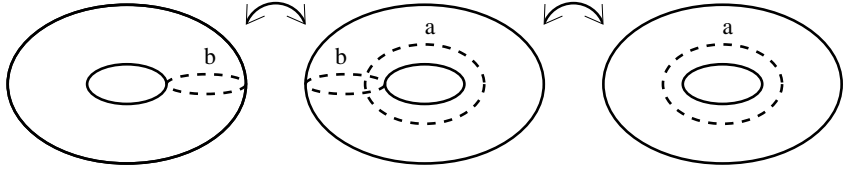}
\end{center}
\caption{Salvetti complex for 3 edges in a line} 
\label{Fig5}
\end{figure}

 In \cite{Gro}, Gromov introduced the notion of a hyperbolic metric space and a word hyperbolic group (one whose Cayley graph is hyperbolic).  He showed that such a metric space $X$ has a well-defined ideal boundary, $\partial X$, and any quasi-isometry between two hyperbolic metric spaces $X$  and $Y$ induces a homeomorphism on the boundary.  In particular,  we can associate to a hyperbolic group $G$ a well-definied boundary, $\partial G$, namely the boundary of (any) Cayley graph of $G$.
 
It is possible to define a boundary for a CAT(0) metric space  $X$ analogously by choosing a basepoint $x_0$ and setting
\[
\partial X =\{\textrm{geodesic rays emanating from $x_0$}\}
\]
with an appropriate topology.
 It is easy to verify that this space is independent of choice of basepoint. 
However, it was shown by P.~Bowers and K.~Ruane \cite{BoRu} that in this setting, quasi-isometries need not extend to homeomorphisms of the boundary.  

In \cite{CK}, C.~Croke and B.~Kleiner  addressed the question of whether a group could act geometrically (i.e.,  properly, cocompactly by isometries) on two CAT(0) spaces $X$ and $Y$ with non-homeomorphic boundaries, i.e., whether a CAT(0) group has a well-definied boundary.  They answered this question in the negative.  Their counterexample was the right-angled Artin group $A_\G$ associated to the graph $\G$ consisting of three edges joined in a line.  
We have already seen one CAT(0) space $X$ on which $\AG$ acts geometrically, namely  the universal cover of the Salvetti complex, $X=\tilde \Sa$.  For this group, the Salvetti complex is made up of three tori, each of which is the quotient of a regular Euclidean cube (Figure \ref{Fig5}). To obtain their counterexample, Croke and Kleiner varied the metric on $\Sa$ by replacing  the cube for the center torus with a parallelogram with acute angle $\alpha$.  They showed that for $\alpha \neq \pi/2$, the boundary of the resulting universal covering space $X_\alpha$ is not homeomorphic to  $\partial X$.  This result was refined by J.~Wilson \cite{Wil} who showed the mind-boggling fact that $\partial X_\alpha \ncong \partial X_\beta$ for any $0 < \alpha < \beta \leq \pi/2$.

\subsection{Applications to robotics}\label{robots}

Right-angled Artin groups also appear in the study of certain problems in robotics.
Abrams and Ghrist  (see eg., \cite{Abr} , \cite{AbGh}, \cite{Ghr}) consider the case of robots moving along tracks on a factory floor.  Since two robots cannot occupy the same space at the same time, the possible positions of the robots is the  configuration space of $n$ distinct points on the graph formed by the tracks.  Abrams and Ghrist study the topology  of configuration spaces of finite graphs.  They noticed that  in some cases, the fundamental groups of these configuration spaces, known as ``graph braid groups", had the structure of a right-angled Artin group.  Subsequently, Crisp and Wiest \cite{CrWi} showed that  all graph braid groups naturally embed in right-angled Artin groups and Sabalka \cite{Sab}  showed that conversely, every right-angled braid group embeds in a graph braid group. In \cite{FaSa},  Farley  and Sabalka determine which of the graph braid groups actually are right-angled Artin groups.

In \cite{GhP}, Ghrist and Petersen consider a different problem in robotics.  Suppose that one has a robot made up of many moving parts.  We can look at the space of all possible configurations of this robot.  Ghrist and Peterson show that, under appropriate conditions, this space can be modeled with a locally CAT(0) cube complex called the ``state complex" for the reconfiguarable system.  A good example to keep in mind is the case of a planar robot made up of a finite collection hexagons attached along edges. A hexagon on the outside of the configuration can rotate at one of its vertices to a new position, giving a new configuration of the robot. The state complex $\mathcal S$ for this system is constructed as follows. The vertices of $\mathcal S$ correspond to all possible configurations of the robot (i.e., arrangements of the hexagons).  Edges of $\mathcal S$ correspond to the rotation of a single hexagon from one position to another, and $k$-cubes correspond to sets of $k$ rotations which are ``independent" of each other (in an appropriately defined sense).  

Following arguments of Crisp and Wiest,  Ghrist and Petersen prove that the fundamental group of the state complex of any such system embeds in the right-angled Artin group associated to the graph whose vertices are the generating moves (such as rotations of a single hexagon) and whose edges correspond to pairs of independent moves.

Aside from giving a nice picture of these robotic systems, there are some potentially useful applications of this work (see eg., \cite{GhL}).  For example, in a CAT(0) cube complex, given any two vertices, there is a ``canonical cube path" connecting them.  ( For a cube complex which is only locally CAT(0), there is one such path in each homotopy class.)  Moreover, given any path between these vertices, there is an algorithm for reducing  the given path to the canonical path.  Viewing the vertices in the state complex as configurations of a robotic system, the canonical cube path can then be used to determine the most efficient ways of  moving from one configuration of the robot to another. 

\medskip

We remark in closing, that there are several other classes of groups closely related to right-angled Artin groups.  These include the graph braid groups described above and the ``mock right-angled Artin groups" recently introduced by  R.~Scott in \cite{Sc}.  These groups are interesting in their own right and suggest possible directions for further study.


\end{document}